\documentclass[12pt,amsfonts,amscd]{amsart}

\usepackage{amsfonts}
\usepackage{color}
\usepackage{mathrsfs}
\newtheorem{Theorem}{Theorem}[section]
\newtheorem{Definition}[Theorem]{Definition}
\newtheorem{Proposition}[Theorem]{Proposition}

\newtheorem{Lemma}[Theorem]{Lemma}
\newtheorem{Corollary}[Theorem]{Corollary}
\theoremstyle{remark}
\newtheorem{Example}[Theorem]{Example}

\def\CC{{\Bbb C}}

\def\NN{{\Bbb N}}
\def\PP{{\Bbb P}}

\def\l{\lambda}

\def\ulim{\mathop{\overline{\rm lim}}}

\def\eps{\varepsilon}

\def\PSH{\operatorname{PSH}}

\def\be{\begin{enumerate}}
\def\ee{\end{enumerate}}
\def\bT{\begin{Theorem}}
\def\eT{\end{Theorem}}
\def\bP{\begin{Proposition}}
\def\eP{\end{Proposition}}
\def\bD{\begin{Definition}}
\def\eD{\end{Definition}}
\def\bE{\begin{Example}}
\def\eE{\end{Example}}
\def\bL{\begin{Lemma}}
\def\eL{\end{Lemma}}
\def\bC{\begin{Corollary}}
\def\eC{\end{Corollary}}

\def\bpp{\begin{proof}}
\def\epp{\end{proof}}
\def\bee{\begin{equation}}
\def\eee{\end{equation}}

\def\Conv{\operatorname{Conv}}
\def\deg{\hbox{\rm deg}\,}
\def\GD{$G_\delta$}
\def\FS{$F_\sigma$}

\def\plp{pluripolar}
\def\psh{plurisubharmonic}
\def\iff{if and only if\ }

\def\beqq{\begin{eqnarray*}}
\def\eeqq{\end{eqnarray*}}
\def\a{\alpha}
\def\b{\beta}
\def\d{\delta}
\def\ve{\varepsilon}

\def\ol{\overline}

\newcommand*\gh[1]{{\hat #1}^G}
\def\span{\text{span}}
\begin{document}
\title[Convergence Sets of Power Series]{On Convergence Sets of Formal Power Series}
\author{Daowei Ma and Tejinder S.~Neelon}

\begin{abstract}
The (projective) convergence set of a divergent formal power series $%
f(x_{1},\dots ,x_{n})$ is defined to be the image in $\PP^{n-1}$ of the set of all $x\in \mathbb{C}^{n}$
such that $f(x_{1}t,\dots ,x_{n}t)$, as a series in $t$, converges
absolutely near $t=0$. We prove that every countable union of closed
complete pluripolar sets in $\PP^{n-1}$ is the convergence set of some
divergent series $f$. The (affine) convergence sets of formal power series
with polynomial coefficients are also studied. The higher-dimensional
results of A.~Sathaye, P.~Lelong, N.~Levenberg and R.E.~Molzon, and of J.~Rib\'on are thus generalized.
\end{abstract}

\address{ dma@math.wichita.edu, Department of Mathematics, Wichita State
University, Wichita, KS 67260-0033, USA}
\address{ neelon@csusm.edu, Department of Mathematics, California State
University, San Marcos, CA 92096-0001, USA}
\maketitle


\section{Introduction}

A formal power series $f(x_{1},\dots ,x_{n})$ with coefficients in $\mathbb{C%
}$ is said to be convergent if it is absolutely convergent in some
neighborhood of the origin in $\mathbb{C}^{n}$. A classical result of Hartogs (see \cite{Ha})
states that a series $f$ converges if
and only if it converges along all directions $\xi \in \mathbb{P}^{n-1}$, {\it 
i.e.}, $f_{\xi }(t):=$ $f($ $\xi _{1}t,\dots ,$ $\xi _{n}t)$ converges, as a
series in $t,$ for all $\xi \in \mathbb{P}^{n-1}$. This can be
interpreted as a formal analog of Hartogs' theorem on separate analyticity.
Since a divergent power series still may converge in certain directions, it
is natural and desirable to consider the set of all such directions. Following Abyankar-Moh
\cite{AM}, we define the convergence set of a \emph{divergent} power series $f$
to be the set of all directions $\xi \in \mathbb{P}^{n-1}$ such that $%
f_{\xi }(t)$ is convergent. For the case $n=2,$ P.~Lelong \cite{Le} proved that
the convergence set of a divergent series $f(x_{1},x_{2})$ is an $F_{\sigma }
$ polar set (i.e. $F_{\sigma }$ set of vanishing logarithmic capacity) in $\mathbb{P}^{1}$, and moreover, every $F_{\sigma }$ polar subset of $\mathbb{P}^{1}$ is contained in the convergence set of a divergent series $%
f(x_{1},x_{2})$. The optimal result was later obtained by
A.~Sathaye (see \cite{Sa}) who showed that the class of convergence sets of divergent  power series $f(x_{1},x_{2})$ is precisely the class of $F_{\sigma }$ polar sets in $\mathbb{P}^{1}$. In this paper we prove that a countable union of closed complete
pluripolar sets in $\mathbb{P}^{n-1}$ is the convergence set of some divergent series. This generalizes the results of P.~Lelong, Levenberg and Molzon, and Sathaye.

We also study convergence sets of power series of the type $f(s,t)=\sum_{j}P_{j }(s)t^{j}$ where the coefficients $P_{j}(s)$ are polynomials with $\deg (P_{j})\leq j$, as in \cite{Ri} and \cite{Pe}. 

Theorems \ref{mainaffine} and \ref{mainproj} are our main theorems, the proofs of which were inspired by \cite{Sa}, and influenced by the methods developed in \cite{Sc3}, \cite{LM}, and \cite{Ri}.

\section{Transfinite Diameter and Capacity}

Let $\mathbb{Z}_{+}$ denote the set of nonnegative integers. Let $n$ be a
positive integer. For $\alpha =(\alpha _{1},...,\alpha _{n})\in \mathbb{Z}%
_{+}^{n}$, let $|\alpha |=\alpha _{1}+...+\alpha _{n}$. Let $\{\alpha
(1),\alpha (2),\dots \}$ be the listing of the elements of $\mathbb{Z}%
_{+}^{n}$ indexed using the lexicographic ordering but with $|\alpha (i)|$
nondecreasing. Set $x^{\alpha }=x_{1}^{\alpha _{1}}\cdots x_{n}^{\alpha
_{n}} $ for $x\in \mathbb{C}^{n}$ and $\alpha \in \mathbb{Z}_{+}^{n}$.

Let $m_{k}={\binom{n+k}{k}}$, the number of monomials of order up to $k$. Let 
$l_{k}=\sum_{q=1}^{k}q(m_{q}-m_{q-1})=n{\binom{n+k}{k-1}}$.

For a finite set $\{s_{1},\dots ,s_{j}\}$ of points in $\mathbb{C}^{n}$, let 
$V(s_{1},\dots ,s_{j})=\det (s_{q}^{\a(p)}))_{1\leq p,q\leq j}$ be the $j$-th
Vandermonde determinant. For a compact set $E$ in $\mathbb{C}^{n}$, let 
\begin{equation*}
V_{j}(E)=\sup \{V(s_{1},\dots ,s_{j}):s_{1},\dots ,s_{j}\in
E\},\;\;d_{k}(E)=(V_{m_{k}}(E))^{1/l_{k}}.
\end{equation*}%
The limit $d(E):=\lim_{k}d_{k}(E)$ exists (\cite{Za}), and is known as the
transfinite diameter of $E$.

Let $\mathscr{P}_k(\CC^n)$ be the set of polynomials on $\CC^n$ of degrees $\le k$.
For a compact set $E\subset \CC^n$ and $p\in \mathscr{P}_k(\CC^n)$, set 
$$|p|_E=\sup\{|p(z)|: z\in E\},$$
$$L_{k,R}(E)=R^{-1}(\sup\{|p|_{\Delta_R}: p\in \mathscr{P}_k(\CC^n), |p|_E\le 1\})^{1/k},$$
and 
$$L_R(E)=\sup_{k} L_{k,R}(E),\;\; c(E)=1/\ulim_{R\to\infty}L_R(E).$$
The quantity $c(E)$ is called the capacity of $E$.

\bL For a compact set $E$ in $\mathbb{C}^{n}$ we have 

(i)$\ c(E)=0$ if and only if $d(E)=0$, and  

(ii) if $n=1$, then $c(E)=d(E)$.
\eL

For a proof of (i), see \cite{LT}. For a proof of (ii), see [Ahlfors 1973, p.~24].

We need to use the following lemma. It appeared, in different forms, in \cite{SW} and \cite{Si}. See, also, \cite{Sc3} and \cite{Za}.

\bL \label{BI} (Bernstein's inequality) Let $E$ be a compact set in $\CC^n$ with $c(E)>0$. Then there is a positive constant $C_E$ such that for every polynomial $p(z)=\sum_{|\a|\le d}a_\a z^\a$, we have $|a_\a|\le C_E^{d}|p|_E$.  \eL  


\section{Some Classes of Pluripolar Sets}

Let $E$ be a Borel subset of $\CC^n$. (Though we do not mention the word ``Borel'' each time, all subsets of $\CC^n$ considered in this paper are assumed to be Borel.) The set $E$ is said to be pluripolar (polar when $n=1$) if for each point $x\in E$ there is a nonconstant plurisubharmonic function $u$ defined in a neighborhood $U$ of $x$ in $\CC^n$ such that $u=-\infty$ on $E\cap U$. The set $E$ is said to be globally pluripolar  if there is a nonconstant plurisubharmonic function $u$ defined on $\CC^n$ such that $E\subset \{y: u(y)=-\infty\}$. Josefson's theorem (answering a question of P.~Lelong) states that $E$ is pluripolar \iff $E$ is globally \plp. 

The set $E$ is said to be complete pluripolar if there is a non-constant plurisubharmonic function $u$ defined on $\CC^n$ such that $E= \{y: u(y)=-\infty\}$. So the set $\{(0, x_2)\in \CC^2: |x_2|<1\}$ and its closure are pluripolar, but not complete \plp. 
A countable union of pluripolar sets is pluripolar. So the set of rationals in the interval $[0,1]$ is polar. It is not complete polar, because each complete \plp\ set is \GD. In $\CC$ each \GD\ polar set is complete polar, which is Deny's theorem (see \cite{De}). 

Following Siciak \cite[P.~2]{Sc3}, we consider families $L,G,H$ of
plurisubharmonic functions:
\begin{eqnarray*}
L(\CC^n) &=&\{u\in \PSH(\mathbb{C}^{n}):\sup_{x\in \mathbb{C}^{n}}(u(x)-\ln
(1+|x|))<\infty \}, \\
G(\CC^n) &=&\exp (L(\CC^n))=\{e^{u}:u\in L(\CC^n)\}, \\
H(\CC^n) &=&\{u\in \PSH(\mathbb{C}^{n}):u\not\equiv 0,u(\lambda x)=|\lambda
|u(x),\forall \lambda \in \mathbb{C},x\in \mathbb{C}^{n}\}.
\end{eqnarray*}

 A \plp\ set $E$ in $\CC^n$ is said to be $L$-complete if there is a non-constant $u\in L(\CC^n)$ such that $E=\{u=-\infty\}$. 
A \plp\ set $F$ in $\CC^n$ is said to be $H$-complete if there is a $w\in H(\CC^n)$ such that $F=\{x: w(x)=0\}$.

It follows from the one-to-one correspondence (see \cite[Prop.~2.7]{Sc3})
\bee\label{eq2}H(\CC\times \CC^n) \ni f(x_0,x)\mapsto  \log f(1,x)\in L(\CC^n)\eee                
between functions of the class $H$ of $n+1$ variables and the functions of the class $L$ of $n$ variables that each 
$H$-complete \plp\ set in $\CC\times\CC^n$ induces a unique $L$-complete \plp\ set in $\CC^n$, and that each $L$-complete \plp\ set in $\CC^n$ is induced by a (not necessarily unique) $H$-complete \plp\ set in $\CC\times\CC^n$.

Let $|x|=(|x_1|^2+\cdots+|x_n|^2)^{1/2}$. Recall that $\mathscr{P}_k(\CC^n)$ is the set of polynomials on $\CC^n$ of degrees $\le k$. Let $\mathscr{H}_k(\CC^n)$ be the set of homogeneous polynomials on $\CC^n$ of degree $k$. 
Let 
$$Q(\CC^n)=\{(p,k):  k\in \NN, p\in\mathscr{P}_k(\CC^n)\},$$
$$|(p(x),k)|=|p(x)|^{1/k},\;\|(p(x),k)\|=|p(x)|^{1/k}/(1+|x|^2)^{1/2},$$
$$|(p,k)|_K=\sup \{|(p(x),k)|: x\in K\},$$
and
$$\|(p,k)\|_K=\sup \{\|(p(x),k)\|: x\in K\},\; \|(p,k)\|=\|(p,k)\|_{\CC^n}.$$
Let 
$$\Gamma(\CC^n)=\{(h,k):  k\in \NN, h\in\mathscr{H}_k(\CC^n)\},$$
$$\|(h(x),k)\|=|h(x)|^{1/k}/|x|,$$
and $$\|(h,k)\|_K=\sup \{\|(h(x),k)\|: x\in K, x\not=0\},\; \|(h,k)\|=\|(h,k)\|_{\CC^n}.$$

\bD Let $F\subset \CC^n$, $F\not=\emptyset$, $x\in\CC^n$, and $0\le r\le 1$. Define
\beqq \tau_H(x, F, r)&=&\inf\{\|(h,k)\|_F: (h,k)\in\Gamma(\CC^n),\;|h(x)|^{1/k}\ge r|x|,\; \|(h,k)\|\le 1\},\\
T_H(x,F)&=&\sup\{r: \tau_H(x, F, r)=0\},\\
\tau_L(x, F, r)&=&\inf\{\|(h,k)\|_F: (h,k)\in Q(\CC^n),\;\|(h(x),k)\| \ge r,\; \|(h,k)\|\le 1\},\\
T_L(x,F)&=&\sup\{r: \tau_L(x, F, r)=0\}.
\eeqq
For the empty set, we define $\tau_H(x, \emptyset, r)=\tau_L(x, \emptyset, r)=0$, and $T_H(x,\emptyset)=T_L(x,\emptyset)=1$.
\eD

It is clear that if $E\subset F$, then $\tau_L(x, E, r)\le \tau_L(x, F, r)$ and $T_L(x,E)\ge T_L(x,F)$.

\bL \label{glo} Let $u\in H(\CC^n)$ be continuous, with $\sup\{u(x): |x|=1\}=1$, and let $F=\{x\in\CC^n: u(x)=0\}$. 
Then for each $x\in \CC^n\setminus F$, $T_H(x,F)\ge u(x)/|x|$. 
\eL

\bpp Fix $x\in \CC^n\setminus F$. Then $x\not=0$, since $u(0)=0$. Let $r$ be a positive number such that $r< u(x)/|x|\le 1$, and let $\d\in (0,r)$. Let $\phi(x)=\max(u(x),\d|x|)$. Then $\phi$ is a continuous function in $H(\CC^n)$. By \cite[Prop.~2.10]{Sc3}, for all $y\in \CC^n$,
$$\phi(y)=\sup\{|h(y)|^{1/k}: (h,k)\in \Gamma(\CC^n),\; |h(z)|^{1/k}\le \phi(z)\;\forall z\in\CC^n\}.$$
Thus there is an $(h,k)\in \Gamma(\CC^n)$ such that $|h(z)|^{1/k}\le \phi(z)\;\forall z\in\CC^n$, and $r|x|<|h(x)|^{1/k}\le \phi(x)$, and therefore
$$|h(x)|^{1/k}> r|x|,\; \|(h,k)\|\le 1,\;\|(h,k)\|_F\le \d.$$
It follows that $\tau_H(x, F, r)\le \d$. Hence $\tau_H(x, F, r)=0$ for each $r<u(x)/|x|$. Therefore, $T_H(x,F)\ge u(x)/|x|$.
\epp

Since (\ref{eq2}) is a one-to-one correspondence, each $L$-complete \plp\ set in $\CC^n$ is related to a  $H$-complete \plp\ set in $\CC\times\CC^n$.

\bP Let $E=\{v=0\}$ be an $L$-complete \plp\ set in $\CC^n$ with $v\in L(\CC^n)$ such that the function $u(x_0,x):=|x_0|\exp(v(x/x_0))$ defined on $\{x_0\not=0\}$ extends to be a continuous function on $\CC\times \CC^n$. 
Then for each $x\in \CC^n\setminus E$, $T_L(x,E)\ge (1+|x|^2)^{-1/2}\exp(v(x))$.
\eP

\bpp This is a consequence of the previous lemma.\epp



\bL \label{tec} Let $E=\{g=0\}$ be a closed $L$-complete \plp\ set in $\CC^n$ with $g\in G(\CC^n)$ such that $\sup\{(1+|y|^2)^{-1/2}g(y): y\in\CC^n\}=1$. Then for each $x\in \CC^n\setminus E$, and each compact set $K$, $T_L(x, E\cap K)\ge (1+|x|^2)^{-1/2}g(x)$. \eL

\bpp  If $E\cap K=\emptyset$, then the desired inequality clearly holds since $T_L(x,\emptyset)=1$. Fix $x\in \CC^n\setminus E$ and a compact set $K$ with $E\cap K\not=\emptyset$. Let $r>0$ be such that $r<(1+|x|^2)^{-1/2}g(x)$. 
Let $\eta$ be a positive number with $\eta<r$. 
Let $\l$ be a positive number that is less than the distance between the closed set $\{y: g(y)\ge \eta \}$ and the compact set $ K\cap E$, and that is so small that
\bee\label{sqrteta}(\l+\eta)^{1-\l}<\sqrt\eta.
\eee
Let
$$\omega(y)=\left\{
\begin{array}{ll}
c_n\exp(-1/(1-|y|^2)),&\mbox{if}\,|y|<1,\\
0,&\mbox{if}\, |y|\ge 1,\end{array}
\right.\;\;\int \omega(y)\,dy=1.$$
For $\mu>0$, let $g_\mu(y)=\int g(y+\mu z)\omega(z)\,dz$. Then $g_\mu\in G(\CC^n)$, $g_\mu$ is $C^\infty$ and positive, and $g_\mu\downarrow g$ as $\mu\downarrow 0$.  If $y\in K\cap E$, and if $|z|<1$, then $y+\l z\not\in\{y: g(y)\ge \eta \}$, and hence $g(y+\l z)<\eta$. It follows that $g_\lambda(y)=\int g(y+\lambda z)\omega(z)\,dz<\int \eta\omega(z)\,dz=\eta$. For each $y\in \CC^n$,
$$g_\lambda(y)\le \int (1+|y+\lambda z|^2)^{1/2}\omega(z)\,dz\le (1+(|y|+\lambda)^2)^{1/2},$$
and hence 
$$g_{\lambda}(y)\le (1+\lambda)(1+|y|^2)^{1/2}.$$
As in  \cite[p.~17]{Sc3}, we define a function $\phi_\lambda\in H(\CC\times\CC^n)$ by
$$\phi_\lambda(y_0,y)=\left\{
\begin{array}{ll}
|y_0|(\lambda+g_\lambda(y/y_0))^{1-\lambda}+\lambda(|y_0|^2+|y|^2)^{1/2}, &\mbox{if}\, y_0\not=0,\\
\lambda|y|,&\mbox{if}\, y_0=0,\end{array}
\right.$$
and define a function $\psi_\lambda\in G(\CC^n)$ by 
$$\psi_\lambda(y)=\phi_\lambda(1,y)=(\lambda+g_\lambda(y))^{1-\lambda}+\lambda(1+|y|^2)^{1/2}.$$
Then $\psi_\lambda$ is $C^\infty$, and $\phi_\lambda$ is continuous.

By \cite[Prop.~2.10]{Sc3}, 
$$\phi_\lambda(y_0, y)=\sup\{|h(y_0,y)|^{1/k}\},$$
where the supremum  is taken over all $(h,k)\in \Gamma(\CC\times\CC^n)$ such that $|h(z_0,z)|^{1/k}\le \phi_\lambda(z_0,z)\;\forall (z_0,z)\in\CC\times\CC^n$. It follows that 
\bee\label{eq5}\psi_\lambda(x)=\sup\{|p(x)|^{1/k}: (p,k)\in Q(\CC^n),\; |p(y)|^{1/k}\le \psi_\lambda(y)\;\forall y\in\CC^n\}.\eee
For all $y\in\CC^n$,  
\begin{eqnarray*}
\psi_\lambda(y)&=&(\lambda+g_\lambda(y))^{1-\lambda}+\lambda(1+|y|^2)^{1/2}\\
&\le & (\lambda+ (1+\lambda)(1+|y|^2)^{1/2})^{1-\lambda}+\lambda(1+|y|^2)^{1/2}\\
&<& (1+3\l)(1+|y|^2)^{1/2}.
\end{eqnarray*}
If $y\in K\cap E$, then
\begin{eqnarray*}
\psi_\lambda(y)&=&(\lambda+g_\lambda(y))^{1-\lambda}+\lambda(1+|y|^2)^{1/2}\\
&\le & (\lambda+ \eta)^{1-\lambda}+\lambda(1+|y|^2)^{1/2}\\
&<& \sqrt\eta+\lambda(1+|y|^2)^{1/2}\\
&\le& (\sqrt\eta+\l)(1+|y|^2)^{1/2}.
\end{eqnarray*}
So $|p(z)|^{1/k}\le \psi_\lambda(z)\,\forall z\in\CC^n$ implies  that 
\bee \label{eqn5}\|(p,k)\|\le 1+3\l \text{ \ \ and \ \ } \|(p,k)\|_{ K\cap E}\le \sqrt\eta+\l.\eee
For sufficiently small $\l$, 
$$(\lambda+g_\lambda(x))^{1-\lambda}+\lambda(1+|x|^2)^{1/2}>(1+3\l)r(1+|x|^2)^{1/2},$$
since as $\l$ approaches 0, the difference of the left side minus the right side tends to $g(x)-r(1+|x|^2)^{1/2}>0$.
It follows that for sufficiently small $\l$, 
\bee\label{eqn8}\psi_\lambda(x)>(1+3\l)r(1+|x|^2)^{1/2}.\eee

By (\ref{eq5}), (\ref{eqn5}) and (\ref{eqn8}), we have $\tau(x, E\cap K,r)\le (1+3\l)^{-1}(\sqrt\eta+\l)$.
Letting $\l\to0$, and then $\eta\to0$, yields that $\tau(x, E\cap K,r)=0$. Since this holds for every $r<g(x)(1+|x|^2)^{-1/2}$, it follows that $T_L(x, E\cap K)\ge g(x)(1+|x|^2)^{-1/2}$.
\epp

\bD A pluripolar set $E$ in $\CC^n$ is said to be $J$-complete if for each $x\in \CC^n\setminus E$, and each compact set $K$, $T_L(x, E\cap K)>0$. 
\eD

Note that the empty set is $J$-complete. Also, it is clear that a $J$-complete \plp\ set has to be closed.

\bP \label{closed} Every closed $L$-complete \plp\ set in $\CC^n$ is $J$-complete.\eP

\bpp This is a consequence of Lemma~\ref{tec}.
\epp

\bP \label{intersection} An intersection of $J$-complete \plp\ sets in $\CC^n$ is $J$-complete. A finite union of $J$-complete \plp\ sets in $\CC^n$ is $J$-complete. \eP

\bpp Let $\{E_\a\}_{\a\in\Lambda}$ be a family of $J$-complete \plp\ sets in $\CC^n$ and let $E=\cap_\a E_\a$. Let $K$ be a compact set in $\CC^n$ and let $x\in \CC^n\setminus E$. Then there is a $\b\in\Lambda$ such that $x\not\in E_\b$. Thus $T_L(x, E\cap K)\ge T_L(x, E_\b\cap K)>0$. Therefore, $E$ is $J$-complete.

Let $F_1, \dots, F_m$ be $J$-complete \plp\ sets in $\CC^n$ and let $F=\cup_{j=1}^m F_j$. Let $K$ be a compact set in $\CC^n$ and let $x\in \CC^n\setminus F$. Choose a number $r$ such that $0<r<\min_j T_L(x, F_j\cap K)$. Then $ \tau_L(x, F_j\cap K,r)=0$ for $j=1, \dots, m$. Let $\ve>0$.  Then there are $(h_j, k_j)\in Q(\CC^n)$, $j=1, \dots, m$, such that 
$$\|(h_j, k_j)\|_{F_j\cap K}<\ve,\;\; \|(h_j(x), k_j)\|\ge r,\;\;  \|(h_j, k_j)\|\le 1.$$
Raising each $h_j$ to a suitable power, we may assume that $k_1=\cdots=k_m=k$. Let $h=\Pi h_j$. Then $(h,mk)\in Q(\CC^n)$, and 
$$\|(h, mk)\|_{F\cap K}<\ve^{1/m},\;\; \|(h(x), mk)\|\ge r,\;\;  \|(h, mk)\|\le 1.$$
Thus $\tau_L(x, F\cap K,r)\le \ve^{1/m}$ for each $\ve>0$. It follows that  $\tau_L(x, F\cap K,r)=0$ and $T_L(x, F\cap K)\ge r>0$. Therefore, $F$ is $J$-complete.

\epp

The following theorem is due to A.~Saddulaev, Since his book that includes the theorem has not been published, we include his proof here. We are grateful to him for sending us the statement and proof of the theorem, and to B.~Fridman for translating an explanation message of A.~Saddulaev from Russian to English.

\bT \label{sad} Every complete \plp\ set in $\CC^n$ is $L$-complete.\eT

\bpp Suppose that $E$ is a complete pluripolar set in $\CC^n$. Let $u$ be a \psh\ function such that $E=\{x:u(x)=-\infty \}$. Choose an increasing sequence $\{M_j\}$ of positive numbers such that $\lim M_j=\infty $ and  $M_j\ge \sup_{|z|\le \exp 2^j} u(z)$. For each $j$, define a function $v_j$ by

\begin{equation*}
v_{j}(x)=\left\{ 
\begin{array}{ll}
\max  (2^{-j}(M_j^{-1}u(x)-1), 2^{-j}\log |x|-1), & \text{ \  if \ } |x|< \exp 2^j; \\ 
2^{-j}\log |x|-1, & \text{ \  if \ } |x|\ge  \exp 2^j.\end{array}
\right. 
\end{equation*}
Since for each $\zeta$ on the boundary of the ball $B(0, \exp 2^j)$, 
$$\limsup_{|x|< \exp 2^j, x\to \zeta} (2^{-j}(M_j^{-1}u(x)-1))\le 0=2^{-j}\log |\zeta|-1,$$
the function $v_j$ is \psh\ on $\CC^n$ by the gluing theorem. On each open set with compact closure, all but a finite number of $v_j$ are non-positive. It follows that the sum $v(x):=\sum_{j=1}^{\infty }v_{j}(x)$ is \psh\ (or identically $-\infty$), since the sequence of the partial sums of the series is eventually non-increasing. It is clear that $v_j(x)\le 2^{-j}\log^+ x$ for each $j$, so that $v(x)\le \log^+x$. Thus $v\in L(\CC^n)$ (or $v$ is identically $-\infty$).

Suppose that $y\in E$. Then $u(y)=-\infty$, and $v_j(y)=2^{-j}\log |y|-1$ for each $j$. Thus $v(y)=-\infty$.

Now suppose that $y\in \CC^n\setminus E$ so that $u(y)>-\infty$. Then $v_j(y)>-\infty$ for each $j$. Since
$$\lim_{j\to\infty} 2^{-j}(M_j^{-1}u(y)-1)=0>-1=\lim_{j\to\infty} (2^{-j}\log |y|-1),$$
it follows that there is a positive integer $m=m(y)$ such that 
$$2^{-j}(M_j^{-1}u(y)-1)>-1/2> (2^{-j}\log |y|-1), \text{ \ for } j>m.$$
Thus $y\in B(0, \exp 2^j)$ and $v_j(y)=2^{-j}(M_j^{-1}u(y)-1)$ for $j>m$, and therefore
\beqq v(y)&=&\sum_{j=1}^m v_j(y)+\sum_{j=m+1}^\infty 2^{-j}(M_j^{-1}u(y)-1)\\
&\ge &\sum_{j=1}^m v_j(y)+\sum_{j=1}^\infty 2^{-j}(-M_1^{-1}|u(y)|-1)\\
&= &\sum_{j=1}^m v_j(y)+(-M_1^{-1}|u(y)|-1)>-\infty.
\eeqq
This implies, in particular, that $v$ is not identically $-\infty$. 
It follows that $v\in L(\CC^n)$ and $E=\{x: v(x)=-\infty\}$.
\epp

\bT \label{cj} Every closed complete \plp set in $\CC^n$ is $J$-complete.\eT

\bpp This is a consequence of Proposition~\ref{closed} and Theorem~\ref{sad}.
\epp 
Let $E$ be a non-empty compact \plp\ set in $\CC^n$. Define the extremal function $\Phi_E: \CC^n\to [0,\infty]$ by $\Phi_E(x)=\sup \{|p(x)|^{1/k}: (p,k)\in Q(\CC^n), |(p,k)|_E\le 1\}$.
The $G$-hull of $E$ is defined to be $\gh{E}:=\{x\in\CC^n: \Phi_E(x)<\infty\}$. The $G$-hull of the empty set is defined to be the empty set. Since $\gh{E}=\cup_k\{x: \Phi_E(x)\le k\}$, it follows that $\gh{E}$ is an \FS\ \plp\ set.

A \plp\ set is said to be $G$-complete if it is the $G$-hull of a compact \plp\ set. A compact complete \plp\ set $K$ in $\CC^n$ is $G$-complete, since $\gh{K}=K$ (see \cite{LM}).



\section{Convergence Sets in Affine Spaces}

Consider a series $f\in \CC[s_1,\dots, s_n][[t]]$ of the form $f(s,t)=\sum_{j=0}^\infty P_j(s)t^j$, where $P_j(s)=P_j(s_1,\dots, s_n)$ are polynomials of $n$ variables.
Define
$$\Conv(f)=\{s\in \CC^n: f(s, t)\; \hbox{\rm converges as a power series in $t$}\}.$$

Let $A, B$ be nonnegative integers with $A>0$. A series $f(s,t)=\sum_j P_j(s)t^j$ is said to be in Class $(A,B)$ if $\deg (P_j)\le Aj+B$. 

It is clear that Class $(1,0)$ is a subset of Class $(A, B)$. Suppose that $E=\Conv(f)$ for some $f$ in Class $(A,B)$. Write $f(s,t)=\sum_j P_j(s)t^j$. Set $g(s, t)=t^Nf(s,t^N)$, where $N=A+B$.
Then $g$ is in Class $(1,0)$ and $\Conv(g)=\Conv(f)$. Therefore, the convergence sets for Class $(A,B)$ are exactly the convergence sets for Class $(1,0)$. 

Suppose that $f(s,t)=\sum_{j=0}^\infty P_j(s)t^j$ is in Class $(1,0)$ and $\Conv(f)=\CC^n$. Then, by Hartogs' classical theorem, $f(s,t)$ converges as a power series in $n+1$ indeterminants $s$ and $t$, {\it i.e.}, $f(s,t)$ converges absolutely for $(s,t)$ in some neighborhood of the origin in $\CC^n\times \CC$.  In this case, we say $f$ is a convergent series. Conversely, if $\Conv(f)\not=\CC^n$, then $f(s,t)$ diverges as a power series in $s$ and $t$, {\it i.e.}, $f(s,t)$ converges absolutely in no neighborhood of the origin in $\CC^{n+1}$. In this case, we say $f$ is a divergent series.

\bD A subset $E$ of $\CC^n$ is said to be a {convergence set} in $\CC^n$ if $E=\Conv(f)$ for some divergent series $f$ of Class $(1,0)$.
\eD

\bT \label{basic} Let $E$ be a convergence set in $\CC^n$. Then $E$ is a countable union of $G$-complete \plp\ sets. Hence $E$ is an \FS\ \plp\ set.\eT

\bpp There is divergent series $f(s,t)=\sum_{j=1}^\infty P_j(s) t^j$ of Class $(1,0)$ such that $E=\Conv(f)$. Put, for $m=1,2,3,\dots$,
\bee \label{eq4a}E_m=\{s\in \CC^n: |s|\le m, |P_j(s)|^{1/j} \le m,\; \text{for}\; j=1,2,\dots\}.\eee
Then $E=\cup E_m$. 

Suppose, if possible, that for some positive integer $m$, $c(E_m)>0$. Then, by Bernstein's inequality (Lemma~\ref{BI}), the coefficients $b_{j\alpha}$ of $P_j(s)=\sum b_{j\alpha } s^\a$ satisfy $|b_{j\a}|\le (C_{E_m}m)^{j}$, where $C_{E_m}$ is a constant depending only on $E_m$. It follows that the series $f(s,t)$ is convergent, contradicting the hypothesis. Therefore each $E_m$ is \plp, and $E$ is an \FS\ \plp\ set.

Fix a non-empty $E_m$ and a point $s\in \gh{E_m}$. Then $\gamma:=\Phi_{E_m}(s)<\infty$. Then $|P_j(s)|^{1/j}\le \gamma m$ for all $j$, and hence $s\in \Conv(f)$. Thus $\gh{E_m}\subset E$ for all $m$.  Therefore, $E=\cup \gh{E_m}$, and $E$ is a countable union of $G$-complete \plp\ sets.
\epp

\bT \label{gcomplete} Every $G$-complete \plp\ set in $\CC^n$ is a convergence set. \eT

\bpp The theorem is proved by following the approach in\\
\cite[Theorem 5.6]{LM}. Let $E$ be a non-empty $G$-complete \plp\ set in $\CC^n$. Then $E=\gh{K}$, where $K$ is a non-empty compact \plp\ set. Let $\mathscr F_K$ be the collection of members $(p,k)\in Q(\CC^n)$ such that $k\ge 1$, $p$ has rational coefficients, and $|(p,k)|_K\le1$. Let $\{(p_j,k_j)\}$ be an enumeration of $\mathscr F_K$. Choose a sequence $\{r_j\}$ of positive integers so that the sequence $\{r_jk_j\}$ is strictly increasing. Let $f(s,t)=\sum_{j=1}^\infty p_j(s)^{r_j}t^{r_jk_j}$. Then $f$ is of Class $(1,0)$.

Suppose $s\in E$. Then $\a:=\Phi_K(s)<\infty$. It follows that $|p_j(s)^{r_j}|\le \a^{r_jk_j}$ for all $j$, and hence $s\in\Conv(f)$. Therefore, $E\subset \Conv(f)$.

We now consider a point $s\not\in E$. Then $\Phi_K(s)=\infty$. For each positive integer $m$ there is a $(p,k)\in Q(\CC^n)$ such that $|(p,k)|_K\le 1$ and $|(p(s),k)|>m$, so there is a $j_m$ such that $|(p_{j_m},k_{j_m})|_K\le 1$ and $|(p_{j_m}(s),k_{j_m})|>m$. It follows that the sequence $\{|(p_j(s)^{r_j},r_jk_j)|\}$ is unbounded, and $s\not\in \Conv(f)$. Therefore, $E= \Conv(f)$.
\epp

\bT \label{main} Let $E$ be a countable union of $J$-complete \plp\ sets in $\CC^n$. Then $E$ is a convergence set.\eT

\bpp The set $E$ can be expressed as $E=\cup E_m$, where $\{E_m\}$ is an ascending sequence of $J$-complete \plp\ sets. For each positive integer $m$, we shall construct a sequence $\{(h_{mk}, q_{mk})\}_{k=1}^\infty$ in $Q(\CC^n)$ such that

(i) $|(h_{mk}, q_{mk})|_{\overline B_m\cap E_m}\le1$,

(ii) $\|(h_{mk}, q_{mk})\|\le m$,

(iii) $\cup_{k=1}^\infty\{x: |(h_{mk}(x), q_{mk})|>m/2\}\supset \CC^n\setminus E_m$,

\noindent where $\ol B_m$ is the closed ball in $\CC^n$ of center 0 and radius $m$.

Fix $m$ and suppose that $y\in \CC^n\setminus E_m$. Then $T_L(x, E_m\cap\ol B_m)>0$. Thus there is a positive number $r<1$ such that
$$\inf\{|(p,v)|_{\overline B_m\cap E_m}: (p,v)\in Q(\CC^n),  |(p(y),v)|\ge r, \|(p,v)\|\le1\}=0.$$
Choose a positive rational number $\beta=a/b<1$, where $a, b$ are positive integers, such that $(r/m)^\beta> 1/2$.  There is a member $(p,v)$ of $Q(\CC^n)$ such that 
$$|(p,v)|_{E_m\cap\ol B_m}<m^{-1/\beta}, \; |(p(y),v)|\ge r, \|(p,v)\|\le1.$$
Let $h_{(y)}(x)=p(x)^a m^{v(b-a)}$, and $q_{(y)}=bv$. Then $(h_{(y)},q_{(y)})\in Q(\CC^n)$, and $|(h_{(y)}(x),q_{(y)})|=|(p(x),v)|^\b m^{1-\b}$. We have, for all $x\in \CC^n$,
$$|(h_{(y)}(x),q_{(y)})|\le (1+|x|^2)^{\beta/2} m^{1-\beta}\le m(1+|x|^2)^{1/2},$$
$$|(h_{(y)},q_{(y)})|_{E_m\cap \ol B_m}<m^{-1}m^{1-\beta}=m^{-\beta}\le 1,$$
and
$$|(h_{(y)}(y),q_{(y)})|\ge r^\beta m^{1-\beta}=(r/m)^\beta m>m/2.$$

Put $U_y:=\{x: |(h_{(y)}(x),q_{(y)})|>m/2\}$. Then $U_y$ is an open  neighborhood of $y$. Since the set $\CC^n\setminus E_m$ is open, the open cover 
$\{U_y: y\in \CC^n\setminus E_m\}$ of $\CC^n\setminus E_m$ contains a countable subcover $\{U_{y_k}: k=1,2,\dots\}$. 
Write $(h_{mk},q_{mk})=(h_{(y_k)},q_{(y_k)})$. Then the sequence $\{(h_{mk}, q_{mk})\}_{k=1}^\infty$ satisfies (i), (ii) and (iii).

Let $\{(P_\nu, q_\nu)\}\}$ be a sequence obtained by arranging $\{(h_{mk}, q_{mk})\}$ as a single sequence. By raising $\{(P_\nu, q_\nu)\}\}$ to suitable powers, we assume that $\{q_\nu\}$ is an increasing sequence. 
Put $f(x,t)=\sum_\nu P_\nu(x)t^{q_\nu}$. Then $f$ is of Class $(1,0)$. We shall show that $E=\Conv(f)$.

Suppose that $x\in E$ and $\nu$ is a positive integer. Then  $x\in \overline B_{m_0}\cap E_{m_0}$ for some positive integer $m_0$. Now $(P_\nu, q_\nu)=(h_{mk}, q_{mk})$ for some $m$, $k$. If $m\ge m_0$, then $|(P_\nu(x), q_\nu)|\le 1$; 
if $m< m_0$, then $|(P_\nu(x), q_\nu)|\le m(1+|x|^2)^{1/2}<m_0(1+|x|^2)^{1/2}$. It follows that $\{|(P_\nu(x), q_\nu)|\}$ is a bounded sequence, and hence $x\in \Conv(f)$.

Now suppose that $x\not\in E$. Then for each positive integer $m$, there is a positive integer $k(m)$ such that
$|(h_{m,k(m)}(x), q_{m,k(m)})|>m/2$. The sequence $\{|(h_{m,k(m)}(x), q_{m,k(m)})|\}_{m=1}^\infty$ is unbounded and is a rearranged subsequence of $\{|(P_\nu(x), q_\nu)|\}$. It follows that $\{|(P_\nu(x), q_\nu)|\}$ is an unbounded sequence, and hence $x\not\in \Conv(f)$. Therefore $E=\Conv(f)$, and $E$ is a convergence set.
\epp

\bT \label{mainaffine} Every countable union of closed complete \plp\ sets in $\CC^n$ is a convergence set.
\eT

\bpp This is a consequence of Theorems~\ref{main} and \ref{cj}.
\epp

\bC Every countable union of proper analytic varieties in $\CC^n$ is a convergence set.
\eC

\bC Every countable set in $\CC^n$ is a convergence set.
\eC

\bC A subset of $\CC$ is a convergence set \iff it is an $F_\sigma$ polar set. 
\eC
\bpp This is because each closed polar set in $\CC$ is a complete polar set.
\epp

\section{Convergence Sets in Projective Spaces}

For a formal power series $f(x_{1},\dots ,x_{n})=f(x)\in \mathbb{C}%
[[x_{1},\dots ,x_{n}]]$ and for $x\in \mathbb{C}^{n}$, let $%
f_{x}(t)=f(x_{1}t,\dots ,x_{n}t)\in \mathbb{C}[[t]]$. Since for $\lambda \in 
\mathbb{C}$, $\lambda \not=0,$ the series $f_{x}$ and $f_{\lambda x}$
converge or diverge together, the convergence set of $f$ ( i.e. the set of $%
x $ for which $f_{x}$ converges) can be identified with a subset of the
projective space $\mathbb{P}^{n-1}$.

For a non-zero member $x$ in $\mathbb{C}^{n}$, $[x]$ denotes its image in $%
\mathbb{P}^{n-1}$. For a subset $E$ of $\mathbb{P}^{n-1}$, put $\tilde{E}%
=\{x\in \mathbb{C}^{n}:[x]\in E\}$.

The (projective)\ convergence set of $f$ is defined to be 
$$\Conv_p(f)=\{[x]\in \PP^{n-1}: f_x\;\mbox{converges}\}.$$

\bD A subset $E$ of $\PP^{n-1}$ is said to be a convergence set in $\PP^{n-1}$ if $E=\Conv_p(f)$ for some divergent series $f(x_1,\dots,x_n)$.
\eD

Let $E$ be a non-empty closed set in $\mathbb{P}^{n-1}$. Define $\Psi _{E}:
\mathbb{P}^{n-1}\rightarrow \lbrack 0,\infty ]$ by $\Psi _{E}([x])=\sup
\{|h(x)|^{1/q}/|x|:(h,q)\in \Gamma (\mathbb{C}^{n}),\Vert (h,q)\Vert _{
\tilde{E}}\leq 1\}$. The $G$-hull of $E$ is $\gh{E}=\{u\in \mathbb{P}%
^{n}:\Psi _{E}(u)<\infty \}$. The $G$-hull of the empty set is defined to be the empty set. If $E$ is non-pluripolar, then $\gh{E}=\mathbb{%
P}^{n-1}$. If $E$ is pluripolar, then $\gh{E}$ is an $F_{\sigma }$\
pluripolar\ set.

Recall that there are no non-constant \psh\ functions on $\PP^{n-1}$. 

\bD A pluripolar\ set $E$ in $\mathbb{P}^{n-1}$ is said to be complete if there
is a function $h\in H(\mathbb{C}^{n})$ such that $E=\{[x]\in \mathbb{P}%
^{n-1}:h(x)=0\}$. A pluripolar set $F$ in $\mathbb{P}^{n-1}$ is
said to be $G$-complete if $F=\gh{E}$ for some closed pluripolar\ set $E$.\eD

The proofs of the following two theorems are very similar to those of Theorems \ref{basic} and \ref{gcomplete}, and hence are omitted.

\bT \label{bas} Let $E$ be a convergence set in $\PP^{n-1}$. Then $E$ is a countable union of $G$-complete \plp\ sets. Hence $E$ is an \FS\ \plp\ set.\eT

\bT Every $G$-complete \plp\ set in $\PP^{n-1}$ is a convergence set. \eT

The set $\Pi $ of all hyperplanes in $\mathbb{P}^{n-1}$ is naturally
isomorphic to $\mathbb{P}^{n-1}$. Each $\Omega \in \Pi $ is isomorphic to $%
\mathbb{P}^{n-2}$, and its complement in $\mathbb{P}^{n-1}$ is isomorphic to 
$\mathbb{C}^{n-1}$. For any two hyperplanes in $\mathbb{P}^{n-1}$, there is
a unitary transformation that maps one to the other.

Fix a positive number $M$. Let 
\beqq  S_1&=&\{[1,0,\dots,0]\},\;\;\text{and for }k=2,\dots,n,\\
S_k&=&\{[x]: |x_1|^2+\cdots+|x_{k-1}|^2\le M^2|x_k|^2, x_{k+1}=\cdots=x_n=0\}.
\eeqq
Put 
\bee \label{km} K_M=\cup_{k=1}^n S_k.\eee Then $\{K_m\}$ is an ascending sequence of closed sets with $\PP^{n-1}=\cup_{m=1}^\infty K_m$.

\bD A subset $E$ of $\PP^{n-1}$ is said to be {\it non-occupying} if there exists $\Omega \in \Pi $ such that $E\cap \Omega =\emptyset$.\eD

\bL \label{noc}If $K$ is a closed non-occupying subset of $\PP^{n-1}$ and if $u\in \PP^{n-1}$, then $K\cup\{u\}$ is non-occupying.\eL

\bpp 
Let $R=\{V\in \Pi: V\cap K=\emptyset\}$ and $S=\{V\in\Pi: u\in V\}$. Then $R$ is a non-empty open set in $\Pi$ and $S$ is a hyperplane in $\Pi$. Thus $R\setminus S$ is non-empty.
\epp

\bL For each $M>0$, the set $K_M$ is non-occupying.       
\eL

\bpp Let $e_1,\dots, e_n$ be the standard basis of $\CC^n$, and let $\eps$ be a sufficiently small positive number.  Let $v_j=e_j+\eps e_{j+1}$ for $j=1,\dots,n-1$.
Put $V_j=\span(v_1,\dots,v_j)$ for $j=1,\dots, n-1$, and $V=V_{n-1}$. Also, let $W_j=\span(e_1,\dots,e_j)$. Note that 
$$ V\cap W_j\subset V_{j-1},\;\;\text{for }j\ge2.$$
Since $S_j\subset W_j$, it follows that $V\cap S_j\subset V_{j-1}$ for $j\ge 2$. It is clear that $V\cap S_1=\emptyset$. For $j\ge 2$ and for sufficiently small $\eps$, since $W_{j-1}\cap S_j=\emptyset$, and since $V_{j-1}$ is close to $W_{j-1}$, we see that $V_{j-1}\cap S_j=\emptyset$. It follows that
$$V\cap K_M=\cup_{j=1}^n (V\cap S_j)\subset \cup_{j=2}^n (V_{j-1}\cap S_j)=\emptyset.$$
Therefore $K_M$ does not intersect the hyperplane $V$.
\epp

\bD A \plp\ set $E$ in $\PP^{n-1}$ is said to be $J$-complete if for each hyperplane $V$, $E\setminus V$ is $J$-complete in $\PP^{n-1}\setminus V$. The set $E$ is said to be globally $J$-complete if for each $[x]\in \PP^{n-1}\setminus E$, $T_H(x, \tilde E)>0$.
\eD

It is clear that each $J$-complete \plp\ set is closed, and that each globally $J$-complete \plp\ set is $J$-complete.

The proof of the following proposition is very similar to that of Proposition~\ref{intersection}, and hence is omitted.

\bP \label{union} An intersection of (globally) $J$-complete \plp\ sets in $\PP^{n-1}$ is (globally) $J$-complete. A finite union of (globally) $J$-complete \plp\ sets in $\PP^{n-1}$ is (globally) $J$-complete.
\eP

\bP Let $E\subset \PP^{n-1}$ be the zero locus of a continuous function $h\in H(\CC^n)$. Then $E$ is a globally $J$-complete \plp\ set in $\PP^{n-1}$.
\eP

\bpp This is a consequence of Lemma~\ref{glo}.
\epp

\bT \label{closed5} Every closed complete \plp\ set in $\PP^{n-1}$ is $J$-complete.\eT
\bpp This is a consequence of Theorem~\ref{cj}.
\epp

\bP \label{glob} A proper algebraic variety in $\PP^{n-1}$ is a global $J$-complete \plp\ set.
\eP

\bpp Let $E$ be a proper algebraic variety in $\PP^{n-1}$. Then there are members $(h_j, q_j)$ of $\Gamma(\CC^n)$, $j=1,\dots, k$, such that
$$E=\{[x]\in \PP^{n-1}: h_1(x)=\cdots=h_k(x)=0\}.$$
Let $h=\sum_{j=1}^k |h_j|^{1/q_j}$. Then $h\in H(\CC^n)$, $h$ is continuous, and $E=\{h=0\}$. By Proposition~\ref{glob}, $E$ is globally $J$-complete.
\epp

For $[x]\in \PP^{n-1}$ and $S\subset \PP^{n-1}$, we define $T_H([x],S)$ to be $T_H(x, \tilde S)$. If $W$ is a hyperplane in $\PP^{n-1}$, and if $z$ and $S$ lie in $\PP^{n-1}\setminus W\cong \CC^{n-1}$, we observe that $T_H(z,S)=0$ \iff $T_L(z,S)=0$.

\bL \label{key} Let $E$ be a $J$-complete \plp\ set in $\PP^{n-1}$, let $K$ be a non-occupying closed set in $\PP^{n-1}$, let $[y]\in \PP^{n-1}\setminus E$, and let $m$ be a real number $\ge 1$. Then there exists an $(h, q)\in \Gamma(\CC^n)$ such that
$$\|(h,q)\|\le m, \;\;\|(h,q)\|_{E\cap K}\le 1,\;\; \|(h(y),q)\|>m/2.$$
\eL

\bpp By Lemma~\ref{noc}, $K\cup\{[y]\}$ is non-occupying, hence there is a hyperplane $V$ such that $\Omega:=\PP^{n-1}\setminus V\supset(K\cup\{[y]\})$. Since $E\cap \Omega $ is a $J$-complete \plp\ set in $\Omega$, we see that $T_L([y], E\cap K)>0$. It follows that $T_H([y], E\cap K)>0$.
Thus there is a positive number $r<1$ such that $\tau _{H}(x,F,r)=0$, {\it \ i.e.},
$$\inf\{\|(p,v)\|_{E\cap K}: (p,v)\in \Gamma(\CC^n),  \|(p(y),v)\|\ge r, \|(p,v)\|\le1\}=0.$$

Choose a positive rational number $\beta=a/b<1$, where $a, b$ are positive integers, such that $(r/m)^\beta> 1/2$.  There is a $(p,v)\in \Gamma(\CC^n)$ such that 
$$\|(p,v)\|_{E\cap K}<m^{-1/\beta}, \; \|(p(y),v)\|\ge r,\; \|(p,v)\|\le1.$$
Let $u=\ol y/|y|$. Then $u$ is a unit vector, and $\langle y,u\rangle:=u_1y_1+\cdots+u_ny_n=|y|$. Put $h(x)=p(x)^a (m\langle x,u\rangle)^{v(b-a)}$, and $q=bv$. Then $(h,q)\in \Gamma(\CC^n)$, and $\|(h(x),q)\|=\|(p(x),v)\|^\b m^{1-\b}$. We have, for all $x\in \CC^n$,
$$\|(h(x),q)\|\le  m^{1-\beta}\le m,$$
$$\|(h,q)\|_{E\cap K}<m^{-1}m^{1-\beta}=m^{-\beta}\le 1,$$
and
$$\|(h(y),q)\|\ge r^\beta m^{1-\beta}=(r/m)^\beta m>m/2.$$
\epp

\bT \label{proj} Let $E$ be a countable union of $J$-complete \plp\ sets in $\PP^{n-1}$. Then $E$ is a convergence set.\eT

\bpp As in the proof of Theorem~\ref{main}, it is enough to construct a sequence $\{(P_{\nu },q_{\nu })\}_{\nu =1}^{\infty }$ in $\Gamma (\mathbb{C}^{n})$,
with $\{q_{\nu }\}$ strictly increasing, such that $[x]\in E$ if and
only if the sequence $\{\Vert (P_{\nu }(x),q_{\nu })\Vert \}$ is bounded,
because then $E$ is the convergence set of $f(x)=\sum_{\nu }P_{\nu }(x)$.

Since, by Proposition~\ref{union}, the union of a finite number of $J$%
-complete pluripolar\ sets is $J$-complete, we can assume that $E=\cup E_{m}$%
, where $\{E_{m}\}$ is an ascending sequence of $J$-complete pluripolar\
sets in $\mathbb{P}^{n-1}$. Recall that $\mathbb{P}^{n-1}=\cup K_{m},$ where 
$K_{m},m=1,2,3,...,$ is the ascending sequence of closed non-occupying sets
in $\mathbb{P}^{n-1}$ defined in (\ref{km}). For each positive integer $m$,
we shall construct a sequence $\{(h_{mk},q_{mk})\}_{k=1}^{\infty }$ in $%
\Gamma (\mathbb{C}^{n})$ such that

(i) $\|(h_{mk}, q_{mk})\|_{K_m\cap E_m}\le1$,

(ii) $\|(h_{mk}, q_{mk})\|\le m$,

(iii) $\cup_{k=1}^\infty\{[x]\in\PP^{n-1}: \|(h_{mk}(x), q_{mk})\|>m/2\}\supset \PP^{n-1}\setminus E_m$.

Fix $m$ and suppose that $[y]\in \PP^{n-1}\setminus E_m$. By Lemma~\ref{key}, there exists an $(h_{[y]}, q_{[y]})\in \Gamma(\CC^n)$ such that
$$\|(h_{[y]},q_{[y]})\|\le m, \;\;\|(h_{[y]},q_{[y]})\|_{E_m\cap K_m}\le 1,\;\; \|(h_{[y]}(y),q_{[y]})\|>m/2.$$

Put $U_{[y]}:=\{[x]: \|(h_{[y]}(x),q_{[y]})\|>m/2\}$. Then $U_{[y]}$ is an open  neighborhood of $[y]$. Since the set $\PP^{n-1}\setminus E_m$ is open, the open cover 
$\{U_{[y]}: [y]\in \PP^{n-1}\setminus E_m\}$ of $\PP^{n-1}\setminus E_m$ contains a countable subcover $\{U_{[y_k]}: k=1,2,\dots\}$. 
Put $(h_{mk},q_{mk})=(h_{[y_k]},q_{[y_k]})$. Then the sequence $\{(h_{mk}, q_{mk})\}_{k=1}^\infty$ satisfies (i), (ii) and (iii).

Let $\{(P_\nu, q_\nu)\}\}$ be a sequence obtained by arranging $\{(h_{mk}, q_{mk})\}$ as a single sequence. By raising $\{(P_\nu, q_\nu)\}\}$ to suitable powers, we assume that $\{q_\nu\}$ is a strictly increasing sequence. Now (i), (ii) and (iii) imply that the sequence $\{\Vert (P_{\nu }(x),q_{\nu })\Vert \}$ is bounded \iff $[x]\in E$.
\epp

\bT \label{mainproj}Every countable union of closed complete \plp\ sets in $\PP^{n-1}$ is a convergence set.
\eT

\bpp This is a consequence of Theorems~\ref{proj} and \ref{closed5}.
\epp

\bC Every countable union of proper algebraic varieties in $\PP^{n-1}$ is a convergence set.
\eC

\bC Every countable set in $\PP^{n-1}$ is a convergence set.
\eC

\bC A subset of $\PP^1$ is a convergence set \iff it is an $F_\sigma$ polar set. 
\eC

\end{document}